\documentclass[11pt, a4paper]{article}
\usepackage{amsmath}
\usepackage{amsfonts}
\usepackage{amssymb}
\usepackage{ulem}

\usepackage{amsthm}
\usepackage{url}
\usepackage{dcolumn}

\usepackage{cite}
\usepackage{fancyhdr}

\usepackage{graphicx, color}
\usepackage{hyperref}
\usepackage{xcolor}

\newtheorem*{rem*}{Remark}

\begin{document}

\title{Rotationally symmetric tilings with convex pentagons \\belonging to 
both the Type 1 and Type 7 families}
\author{ Teruhisa SUGIMOTO$^{ 1), 2)}$ }
\date{}
\maketitle

{\footnotesize

\begin{center}
$^{1)}$ The Interdisciplinary Institute of Science, Technology and Art

$^{2)}$ Japan Tessellation Design Association

E-mail: ismsugi@gmail.com
\end{center}

}

{\small
\begin{abstract}
\noindent
Rotationally symmetric tilings by a convex pentagonal 
tile belonging to both the \mbox{Type 1} and Type 7 families are introduced. Among 
them are spiral tilings with two-fold and four-fold rotational symmetry. Those 
rotationally symmetric tilings are connected edge-to-edge and have no axis 
of reflection symmetry.
\end{abstract}
}



\section{Introduction}
\label{section1}

\begin{figure}[t]
 \centering\includegraphics[width=12cm,clip]{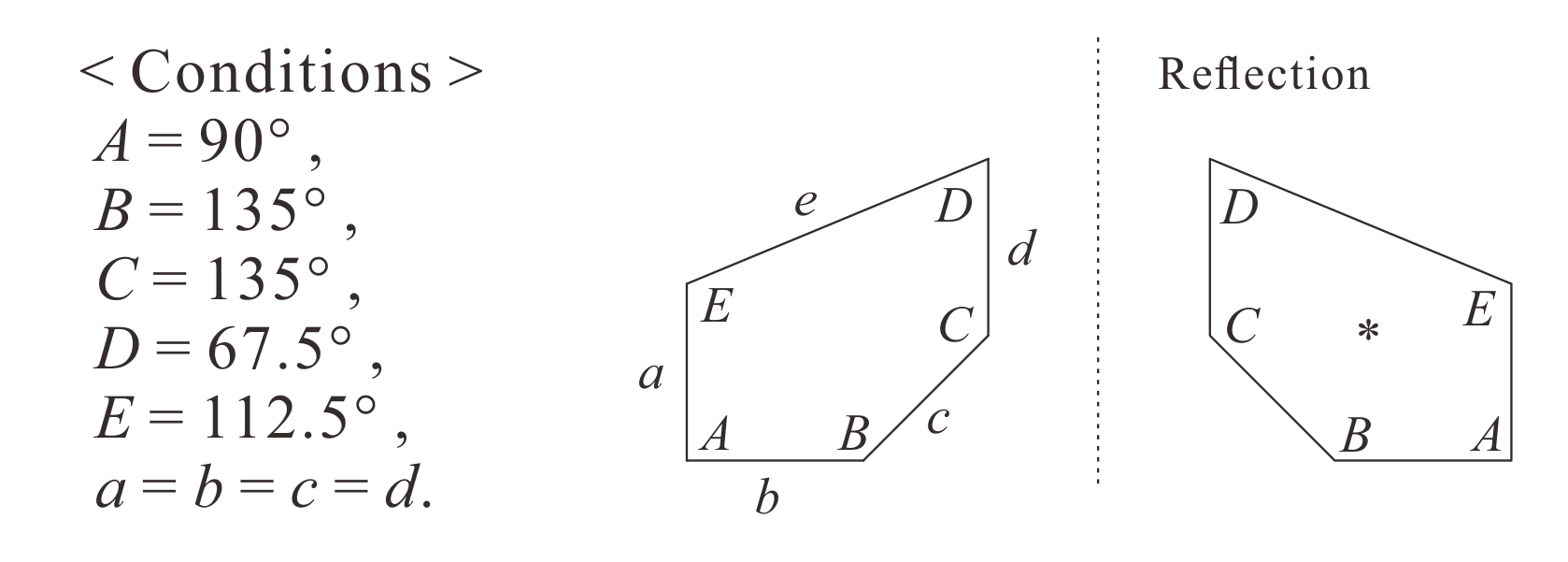} 
  \caption{{\small 
Convex pentagonal tile \mbox{$P(T1 \cap T7)$} that belongs 
to both the Type 1 and Type 7 families
} 
\label{fig01}
}
\end{figure}

To date, fifteen families of convex pentagonal tiles\footnote{
A \textit{tiling} 
(or \textit{tessellation}) of the plane is a collection of sets that are called tiles, 
which covers a plane without gaps and overlaps, except for the boundaries 
of the tiles. The term ``tile'' refers to a topological disk, whose boundary is a 
simple closed curve. If all the tiles in a tiling are of the same size and shape, 
then the tiling is \textit{monohedral}~\cite{G_and_S_1987, Sugimoto_2012, 
Sugimoto_NoteTP, wiki_pentagon_tiling}. In this study, a polygon that admits 
a monohedral tiling is called a \textit{polygonal tile}~\cite{Sugimoto_2012, 
Sugimoto_NoteTP}. Note that, in monohedral tiling, it admits the use of reflected 
tiles.
}$^{,}$\footnote{ 
In May 2017, Micha\"{e}l Rao declared that the complete list of Types of convex 
pentagonal tiles had been obtained (i.e., they have only the known 15 families), 
but it does not seem to be fixed as of March 2020~\cite{wiki_pentagon_tiling}.
}, each of them referred to as a ``Type,'' are known. For example, if the sum 
of three consecutive angles in a convex pentagonal tile is $360^ \circ $, the 
pentagonal tile belongs to the Type 1 family. Convex pentagonal tiles belonging 
to some families also exist. Known convex pentagonal tiles can form periodic 
tiling~\cite{G_and_S_1987, Sugimoto_2012, Sugimoto_NoteTP, wiki_pentagon_tiling}. 
In \cite{Sugimoto_2020_1}, we introduced rotationally symmetric tilings and 
rotationally symmetric tiling-like patterns with a regular convex polygonal 
hole at the center, using convex pentagonal tiles. Then, it showed that the 
convex pentagonal tile that belongs to both the Type 1 and \mbox{Type 7} families 
(see Figure~\ref{fig01}) can form a rotationally symmetric tiling-like pattern 
with a regular convex octagonal hole at the center, and spiral tilings with 
two-fold rotational symmetry. Note that the tiling-like pattern is not 
considered tiling due to the presence of a gap, but is simply called tiling 
in this study. Those rotationally symmetric tilings are connected in an 
edge-to-edge\footnote{ 
A tiling by convex polygons is 
\textit{edge-to-edge} if any two convex polygons in a tiling are either disjoint 
or share one vertex or an entire edge in common. Then other case is 
\textit{non-edge-to-edge}~\cite{G_and_S_1987, Sugimoto_2012, Sugimoto_NoteTP }.
} manner and have no axis of reflection symmetry\footnote{ 
Hereafter, a figure with $n$-fold rotational symmetry without reflection is 
described as $C_{n}$ symmetry. ``$C_{n}$'' is based on the Schoenflies notation 
for symmetry in a two-dimensional point group~\cite{wiki_point_group, 
wiki_schoenflies_notation}. 
}.

\begin{figure}[htb]
 \centering\includegraphics[width=15cm,clip]{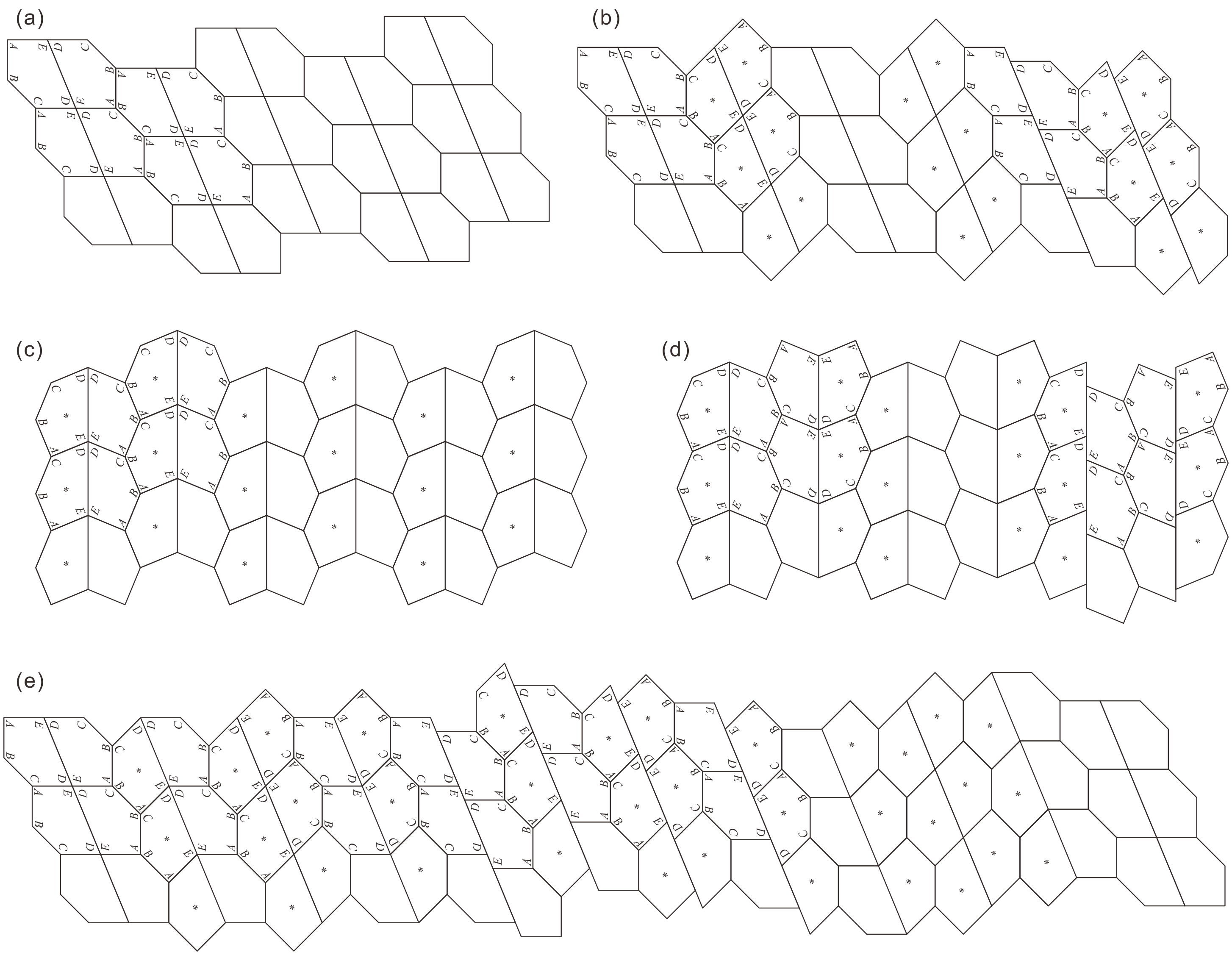} 
  \caption{{\small 
Examples of variations of Type 1 tilings by a convex pentagonal 
tile \mbox{$P(T1 \cap T7)$}
} 
\label{fig02}
}
\end{figure}

\begin{figure}[!h]
 \centering\includegraphics[width=15cm,clip]{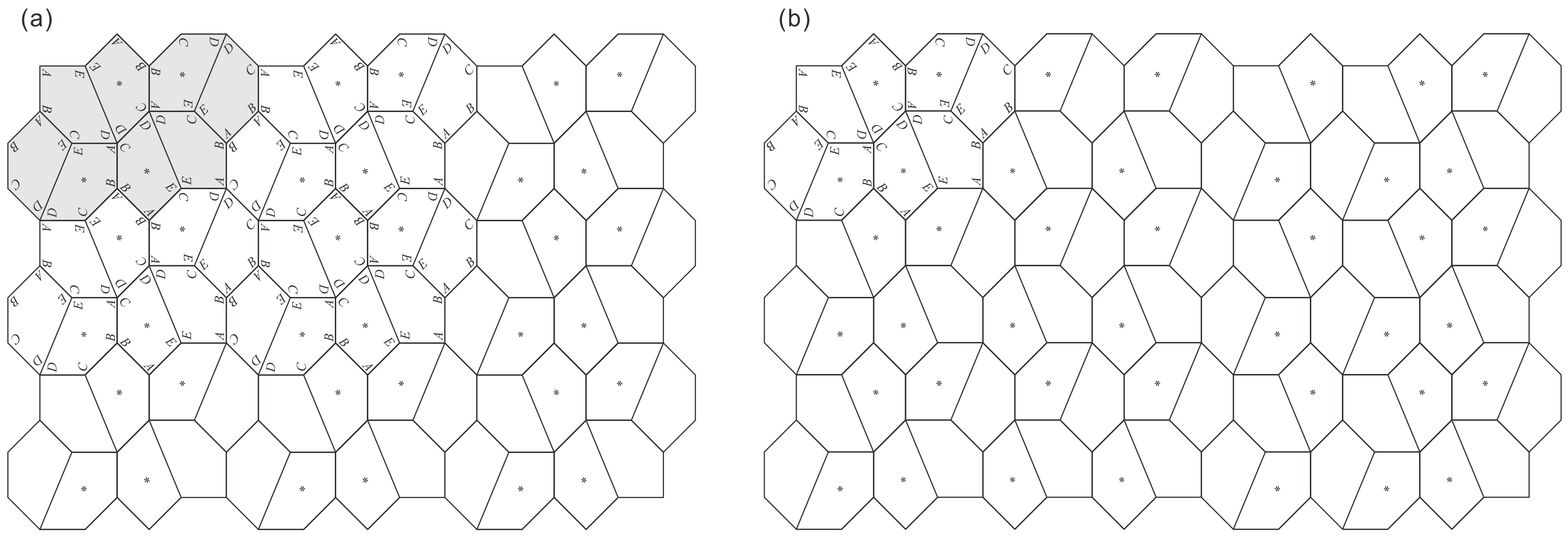} 
  \caption{{\small 
Representative tiling of Type 7 and example of other tiling by 
a convex pentagonal tile \mbox{$P(T1 \cap T7)$}
} 
\label{fig03}
}
\end{figure}

The convex pentagonal tile that belongs to both the Type 1 and Type 7 families 
is the only convex pentagon that satisfies the conditions ``$A = 90^ \circ ,\;
B = C = 135^ \circ ,\;D = 67.5^ \circ ,\;E = 112.5^ \circ ,\;a = b = c = d$'' 
shown in Figure~\ref{fig01}~\cite{Sugimoto_2012, Sugimoto_NoteTP, Sugimoto_2016}. 
Hereafter, let \mbox{$P(T1 \cap T7)$} be the convex pentagonal tile that belongs 
to both the Type 1 and Type 7 families. \mbox{$P(T1 \cap T7)$} can generate 
the representative tiling of Type 1 (see Figure~\ref{fig02}(a)), variations of Type 1 
tilings (i.e., tilings whose vertices are formed only by the relations of 
$A + B + C = 360^ \circ $ and $D + E = 180^ \circ )$ as shown in 
Figures~\ref{fig02}(b)--\ref{fig02}(e), and the representative tiling of Type 7 
(see Figure~\ref{fig03}(a)). It was previously recognized that \mbox{$P(T1 \cap T7)$} 
could generate other tilings as shown in Figure~\ref{fig03}(b)~\cite{Sugimoto_2012}. 
From the research in \cite{Sugimoto_2020_1}, we noticed that it is possible 
to generate more tilings with \mbox{$P(T1 \cap T7)$}. In this study, some of the 
tilings generated by \mbox{$P(T1 \cap T7)$} are introduced.

\renewcommand{\figurename}{{\small Figure.}}

\renewcommand{\figurename}{{\small Figure.}}

\renewcommand{\figurename}{{\small Figure.}}

\renewcommand{\figurename}{{\small Figure.}}
\begin{figure}[htbp]
 \centering\includegraphics[width=10cm,clip]{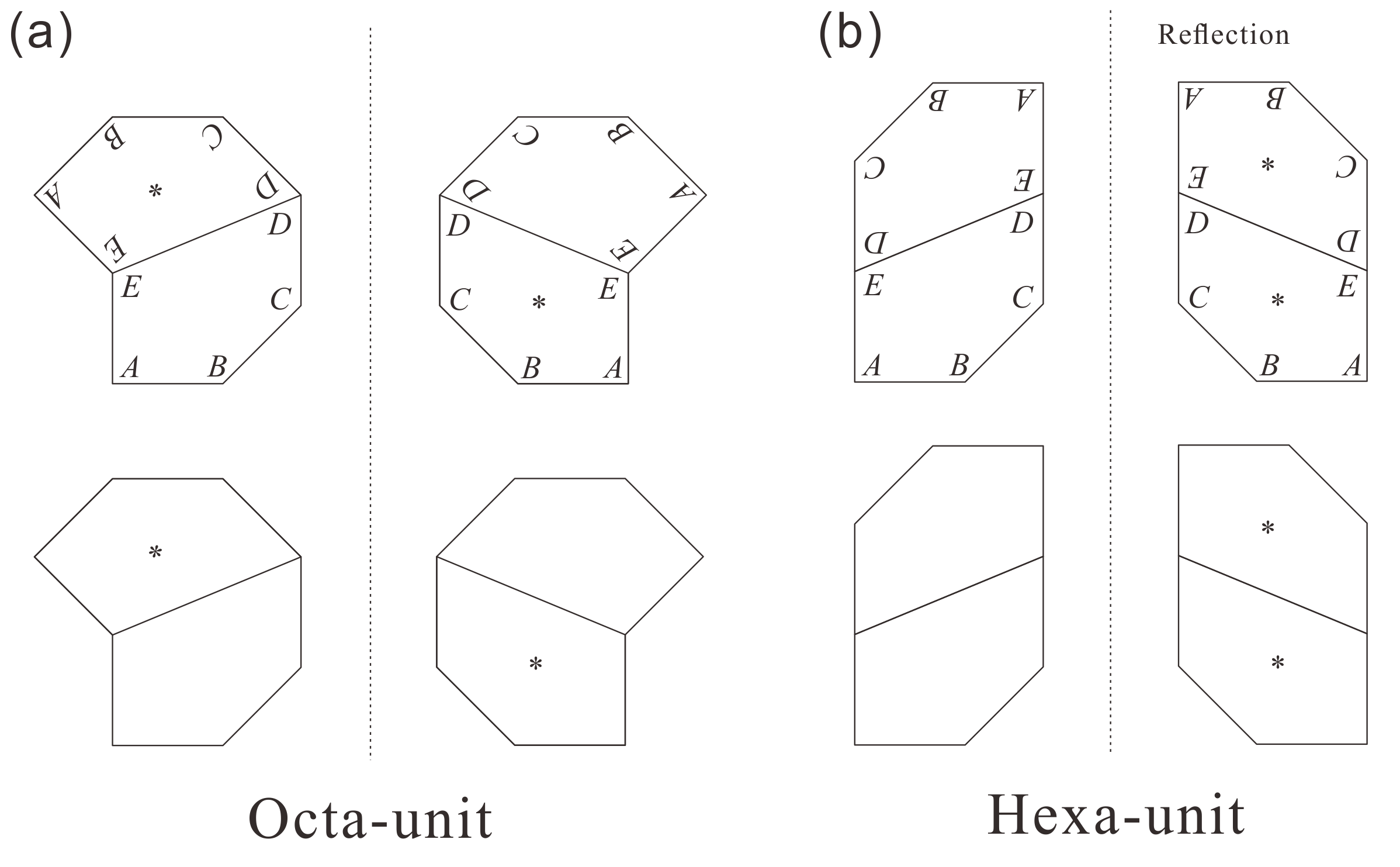} 
  \caption{{\small 
Octa-unit and Hexa-unit formed by two pieces of  
\mbox{$P(T1 \cap T7)$}
} 
\label{fig04}
}
\end{figure}


\section{Octa-unit and Hexa-unit}
\label{section2}

\mbox{$P(T1 \cap T7)$} can also form non-edge-to-edge tilings as shown in 
Figures~\ref{fig02}(b), \ref{fig02}(d), and \ref{fig02}(e) by using the relation 
$D + E = 180^ \circ $, but the tilings by \mbox{$P(T1 \cap T7)$} introduced in this 
study are edge-to-edge. The edge $e$ of \mbox{$P(T1 \cap T7)$} is the only 
edge of different length. As shown in Figure~\ref{fig04}(a), an equilateral 
concave octagon formed by two convex pentagons, connected through 
a line symmetry whose axis is edge $e$, is referred to as the \textit{Octa-unit}. 
As shown in Figure~\ref{fig04}(b), a convex hexagon formed by two convex 
pentagons, connected with rotational symmetry on edge $e$, is referred 
to as the \textit{Hexa-unit}. 

When variations of Type 1 tilings are generated by \mbox{$P(T1 \cap T7)$}, as shown 
in Figure~\ref{fig02}, \mbox{$P(T1 \cap T7)$} can generate tilings with only Octa-units, 
tilings with only Hexa-units, or tilings with both, Octa-units and Hexa-units. The 
representative tiling of Type 7 shown in Figure~\ref{fig03}(a) and other tiling shown 
in Figure~\ref{fig03}(b) are tilings with only Octa-units. 

Note that the candidates of concentration relations that will be used in 
tilings with \mbox{$P(T1 \cap T7)$} are as follows.

\bigskip
\noindent
$D + E = 180^ \circ$, $A + B + C = 360^ \circ$, $2B + A = 360^ \circ$, 
$2C + A = 360^ \circ$, $2E + B = 360^ \circ$, $2E + C = 360^ \circ$, 
$4A = 360^ \circ$, $2A + D + E = 360^ \circ$, $2D + A + B = 360^ \circ$, 
$2D + A + C = 360^ \circ$, $2D + 2E = 360^ \circ$, $4D + A = 360^ \circ$.

\bigskip

\noindent
Edge-to-edge tilings with ``$2E + B = 360^ \circ ,\;2E + C = 360^ \circ 
,\;2D + A + B = 360^ \circ ,\;2D + A + C = 360^ \circ ,\;2D + E = 360^ \circ 
$, or $4D + A = 360^ \circ $'' require the use of Octa-units, and 
edge-to-edge tilings with ``$2A + D + E = 360^ \circ $'' require the use of 
Hexa-units.


\section{Rotationally symmetric tilings}
\label{section3}

\begin{figure}[!h]
 \centering\includegraphics[width=13.5cm,clip]{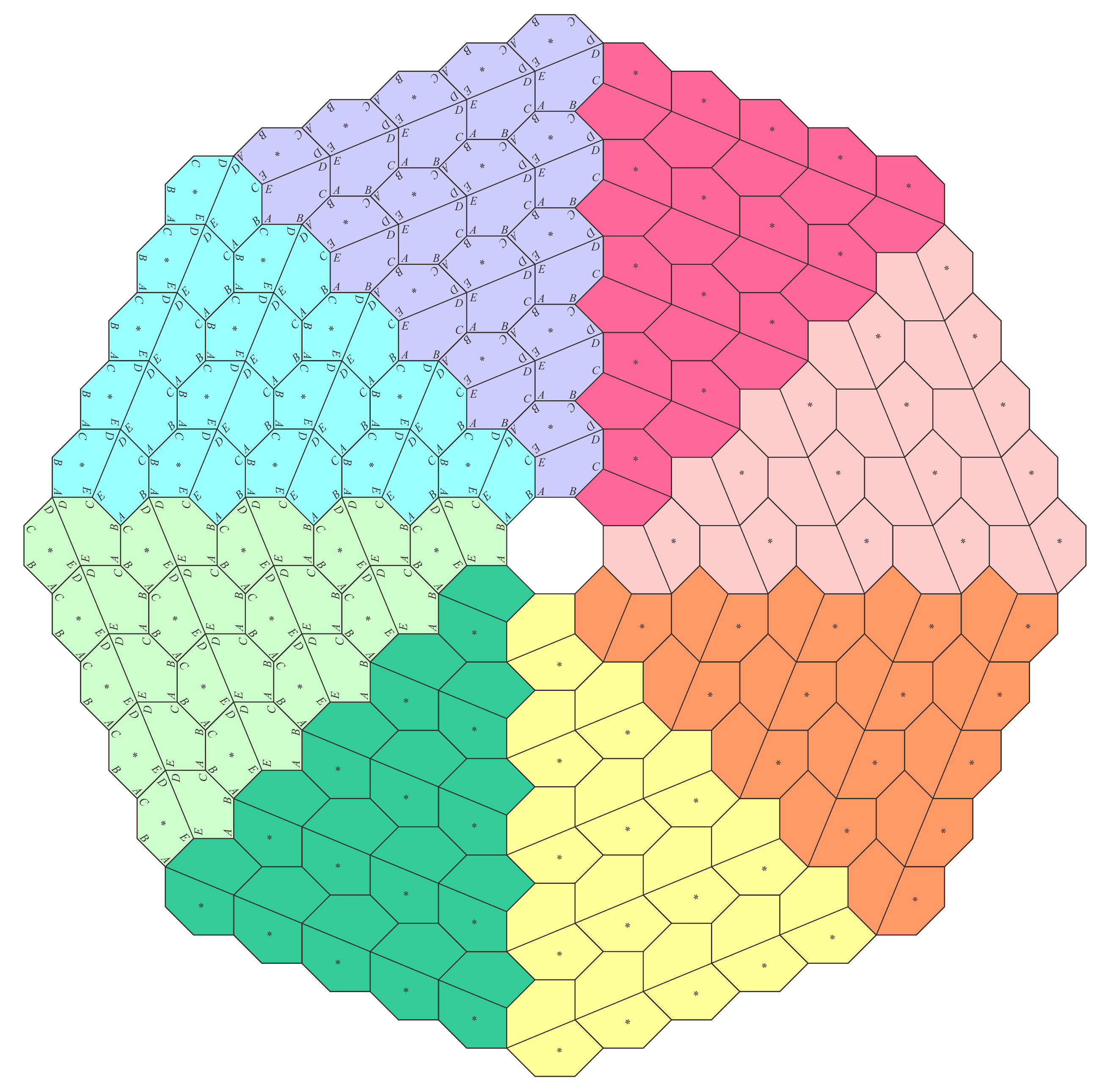} 
  \caption{{\small 
Eight-fold rotationally symmetric tiling with a regular 
convex octagonal hole at the center by a convex pentagonal tile 
\mbox{$P(T1 \cap T7)$} 
} 
\label{fig05}
}
\end{figure}

\begin{figure}[htbp]
 \centering\includegraphics[width=15.5cm,clip]{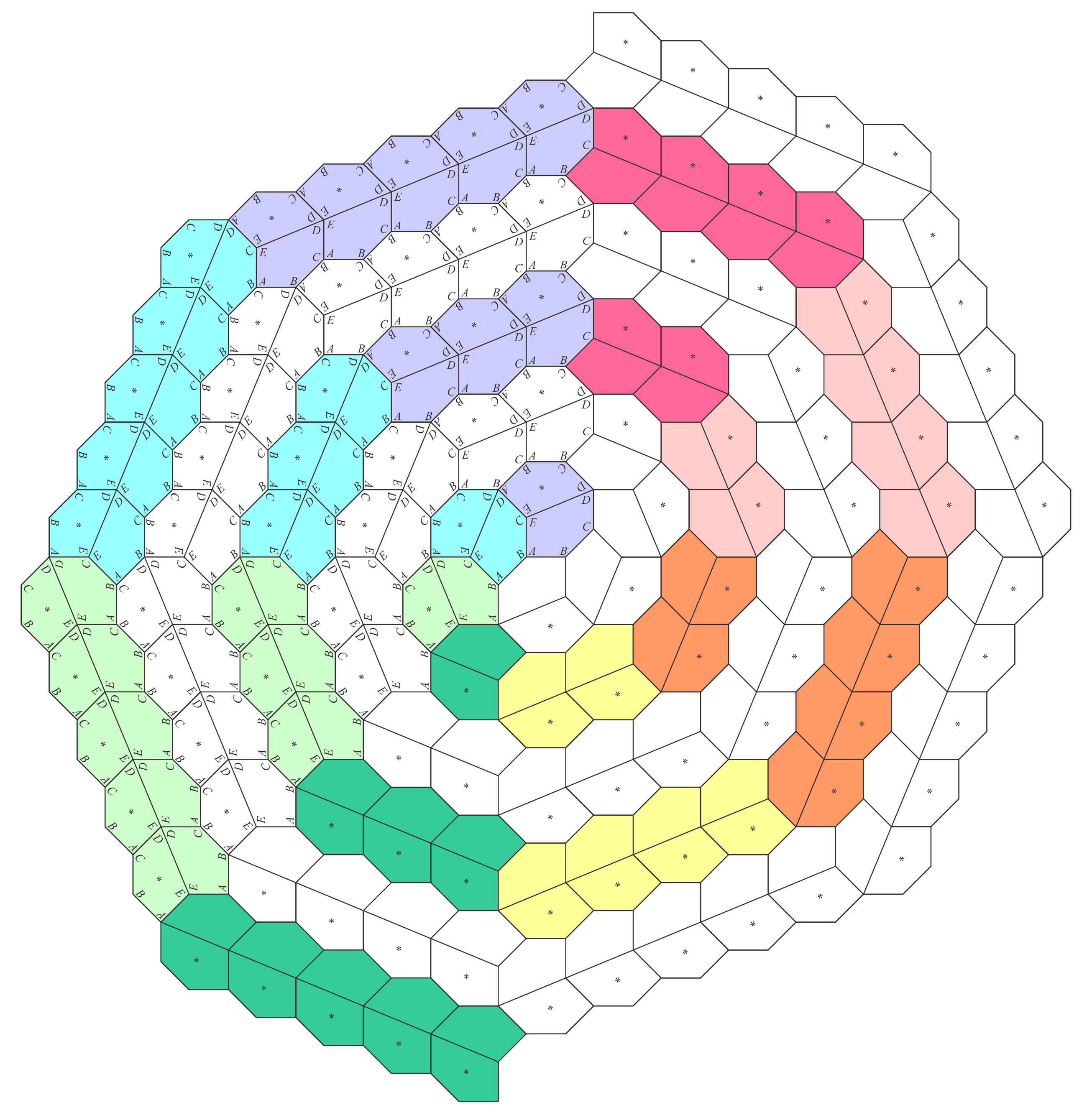} 
  \caption{{\small 
Spiral tiling with two-fold rotational symmetry by a convex pentagonal tile 
\mbox{$P(T1 \cap T7)$} 
} 
\label{fig06}
}
\end{figure}

First, let us introduce rotationally symmetric tilings with only Octa-units. 
In Figure~\ref{fig05}, a rotationally symmetric tiling with $C_{8}$ symmetry, 
with a regular octagonal hole at the center, is depicted as shown in \cite{Sugimoto_2020_1}. 
In Figures~\ref{fig06} and \ref{fig07}, respectively, a spiral tiling with $C_{2}$ 
symmetry and tilings that maintain the spiral structure and extend in one 
direction, are depicted as shown in \cite{Sugimoto_2020_1}. Although a 
tiling is not rotationally symmetric, it is also possible to remove one spiral 
structure and to extend the belts wherein the Octa-units are arranged 
(see Figure 32 in \cite{Sugimoto_2020_1}). 

As mentioned above, \mbox{$P(T1 \cap T7)$} can use the concentration 
$4A = 360^ \circ$ for tiling. When the concentration $4A = 360^ \circ $ 
is formed only by an Octa-unit, there are two types of concave 16-gons 
(hexadecagon) as shown in Figures~\ref{fig08}(a) and \ref{fig08}(b), which 
are reflection images. Hereinafter, the concave 16-gon formed of four 
Octa-units, as shown in Figure~\ref{fig08}, will be referred to as a 
\textit{Hexadeca-4oc-unit}. Hexadeca-4oc-units have four-fold rotational 
symmetry and four axes of reflection symmetry passing through the 
center of the rotational symmetry (hereafter, this property is described 
as $D_{4}$ symmetry\footnote{
``$D_{4}$'' is based on the Schoenflies notation for symmetry in a 
two-dimensional point group~\cite{wiki_point_group, wiki_schoenflies_notation}. 
``$D_{n}$'' represents an $n$-fold rotation axis with $n$ reflection symmetry axes. }). 

David Dailey presented tilings with ``Dented octagons'' on his site~\cite{Dailey_site}. 
His dented octagons correspond to the Octa-unit in Figure~\ref{fig04}(a). The 
figures that he presented are those of very interesting tilings; however, they 
also include patterns that are not expandable without leaving gaps. From his 
results, we can form a spiral tiling with $C_{4}$ symmetry, the center of which is 
the Hexadeca-4oc-unit by \mbox{$P(T1 \cap T7)$}, as shown in Figure~\ref{fig09}. 
(The other tilings with \mbox{$P(T1 \cap T7)$} that can be obtained from his results 
are introduced in the next section.) In addition, we found that \mbox{$P(T1 \cap T7)$} 
can form a spiral tiling with $C_{4}$ symmetry, the center of which is the 
Hexadeca-4oc-unit, as shown in Figure~\ref{fig10}. The tilings presented in 
Figures~\ref{fig09} and \ref{fig10} have the relation of a spiral with reverse 
rotation, and the central Hexadeca-4oc-unit is the same (corresponding to 
Figure~\ref{fig08}(a)). As mentioned above, the Hexadeca-4oc-units have 
$D_{4}$ symmetry; therefore, they can be reversed freely. Considering this 
property, we can determine that when the central Hexadeca-4oc-unit in 
Figure~\ref{fig10} is replaced with that in Figure~\ref{fig08}(b) and the entire 
tiling is reversed, the tiling matches that shown in Figure~\ref{fig09}.

In Figures~\ref{fig11} and \ref{fig12}, non-spiral rotationally symmetric tilings with 
$C_{4}$ symmetry with centers that are Hexadeca-4oc-units, are illustrated. 
The tilings depicted in Figures~\ref{fig11} and \ref{fig12} have a relationship 
in which the central Hexadeca-4oc-unit is reversed.

The rotationally symmetric tilings introduced at the end are tilings with 
$C_{2}$ symmetry as shown in Figures~\ref{fig13} and \ref{fig14} in which an 
Octa-unit and a Hexa-unit are mixed. The center of the tiling depicted in 
Figure~\ref{fig13} is the concentration $4A = 360^ \circ $ that was formed 
with four Hexa-units, and the tiling depicted in Figure~\ref{fig14} has the two 
concentrations $4A = 360^ \circ $ that were formed with two Octa-units 
and two Hexa-units.


\section{Other tilings}
\label{section4}

Using the concentration $4A = 360^ \circ $ formed by mixing an Octa-unit and 
a Hexa-unit, we found the tiling depicted in Figure~\ref{fig15}(a). We also found a 
tiling, as shown in Figure~\ref{fig15}(b), in which Hexadeca-4oc-units are arranged 
side by side in one row and the top and bottom are filled with Octa-units. 

Collections of cases similar to those in Dailey's results, where tilings would not 
be broken by expansion, are shown in Figure~\ref{fig16}~\cite{Dailey_site}. 
In Figure~\ref{fig16}(c), a fusion of the tiling methods shown in Figures~\ref{fig16}(a) 
and \ref{fig16}(b) is depicted, showing that similar tiling can be extended vertically 
and horizontally.


\section{Conclusions}

In this study, we mainly introduced rotationally symmetric tilings that would not 
have been known but for using the convex pentagonal tile \mbox{$P(T1 \cap T7)$} 
belonging to both the Type 1 and Type 7 families. We think that some of 
them can be properly called spiral tiling. In this study, we only introduce 
those tilings with \mbox{$P(T1 \cap T7)$}, and do not fully consider the properties of 
the convex pentagon and other such tilings. (At present, there are many 
aspects that we have not been able to consider; hence, in this article, we 
focused on an introduction of the tiling that we discovered.)

\begin{figure}[htbp]
 \centering\includegraphics[width=15.5cm,clip]{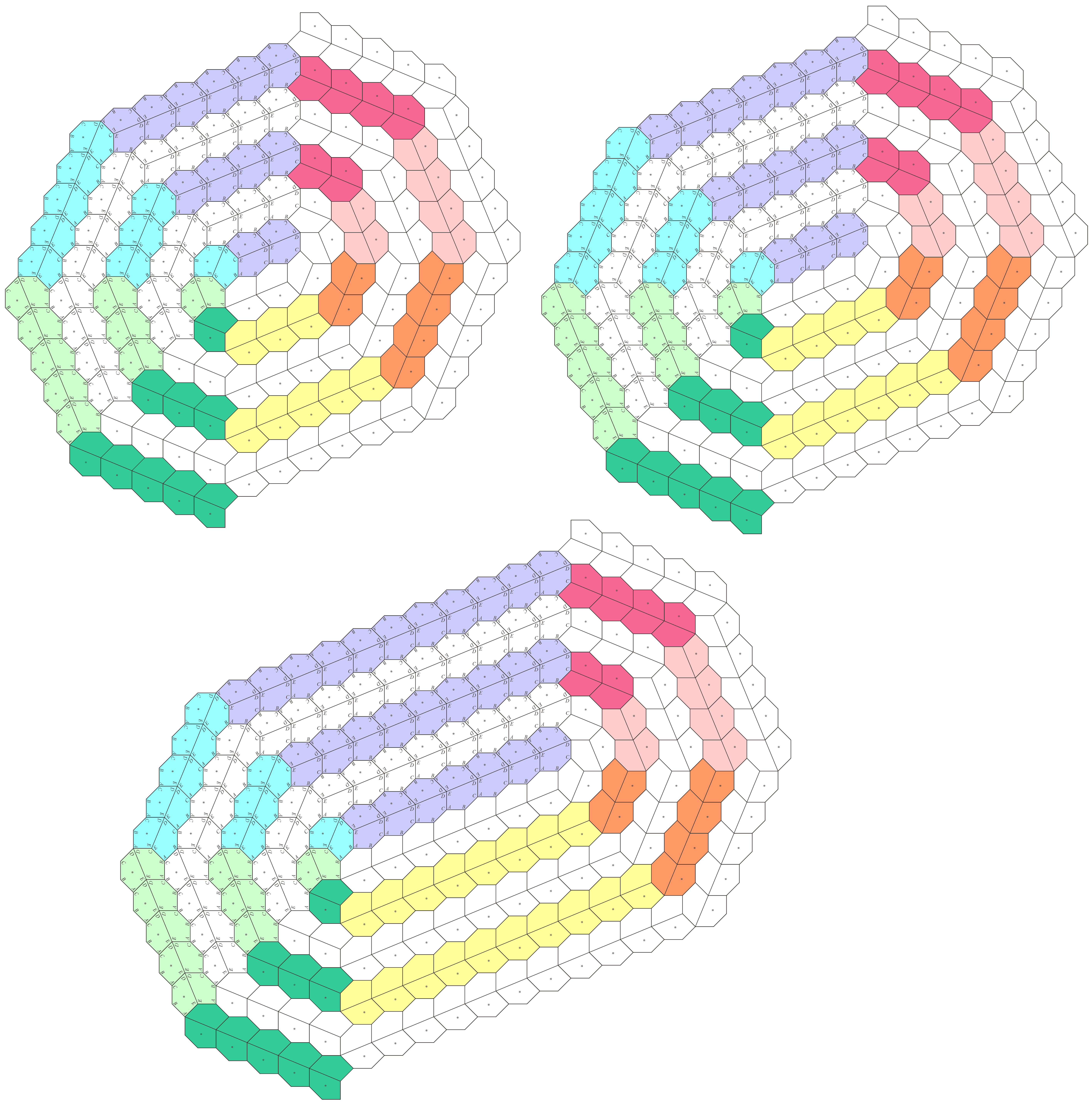} 
  \caption{{\small 
Spiral tilings that maintain the spiral structure and extend in one 
direction by a convex pentagonal tile \mbox{$P(T1 \cap T7)$} 
} 
\label{fig07}
}
\end{figure}

\begin{figure}[htbp]
 \centering\includegraphics[width=12cm,clip]{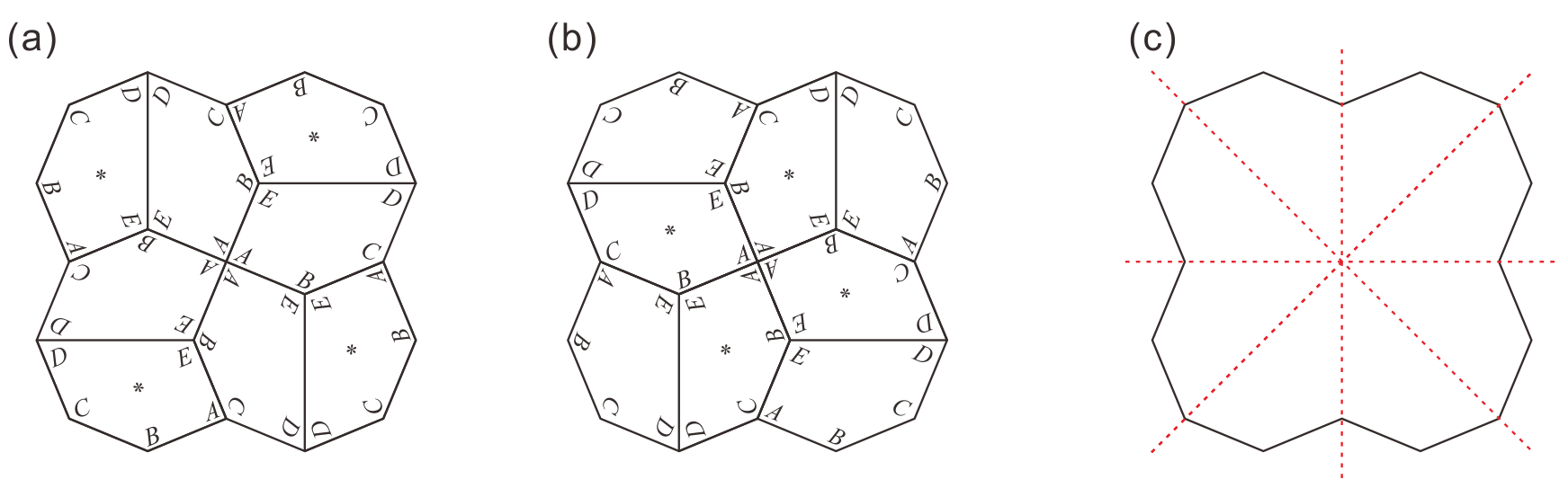} 
  \caption{{\small 
Concave 16-gons  (Hexadeca-4oc-units) that formed of four Octa-units
} 
\label{fig08}
}
\end{figure}

\begin{figure}[htbp]
 \centering\includegraphics[width=15.5cm,clip]{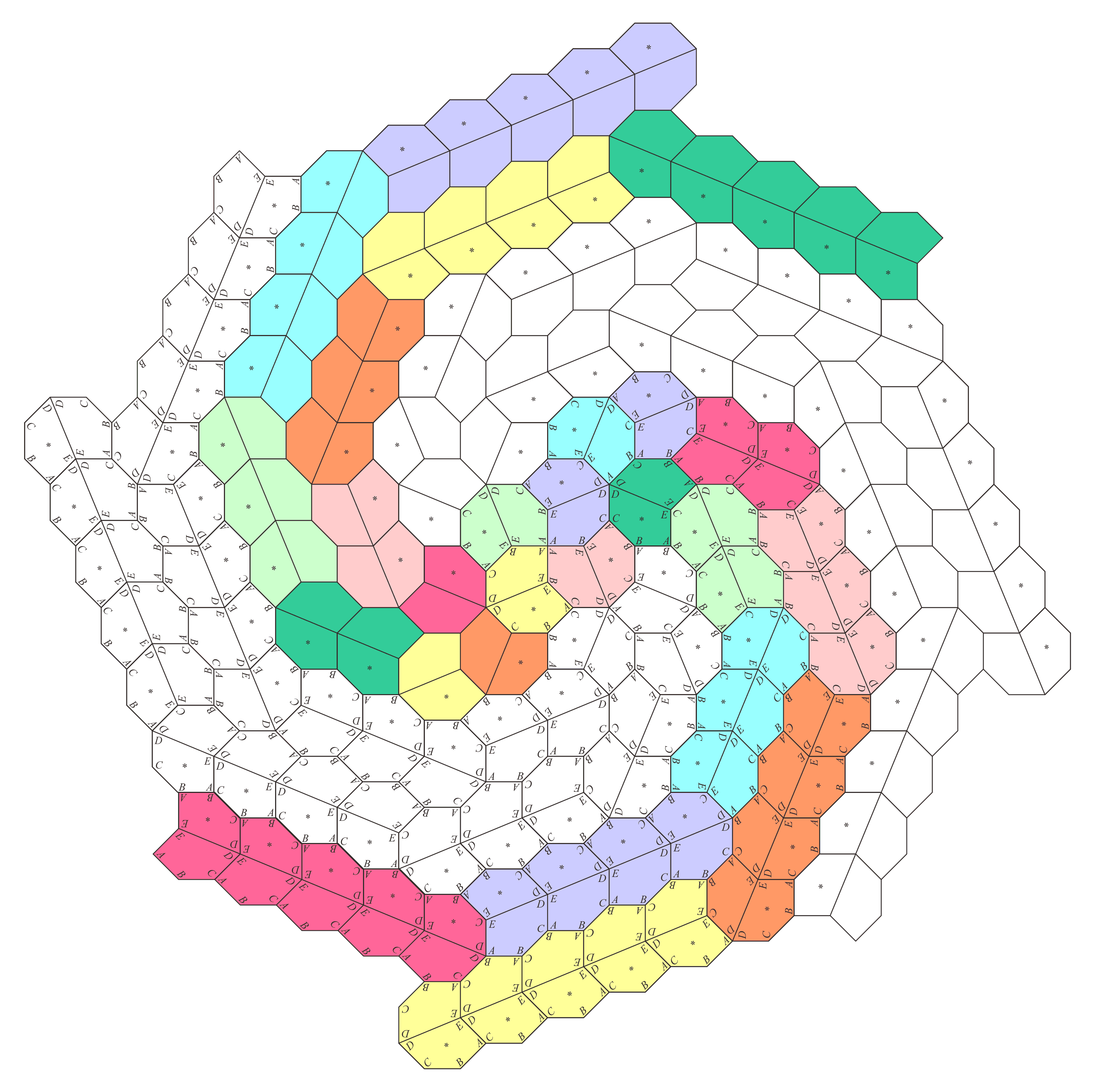} 
  \caption{{\small  
Spiral tiling with four-fold rotational symmetry by a convex pentagonal tile 
\mbox{$P(T1 \cap T7)$} 
based on Dailey's tiling
} 
\label{fig09}
}
\end{figure}

\begin{figure}[!h]
 \centering\includegraphics[width=15.5cm,clip]{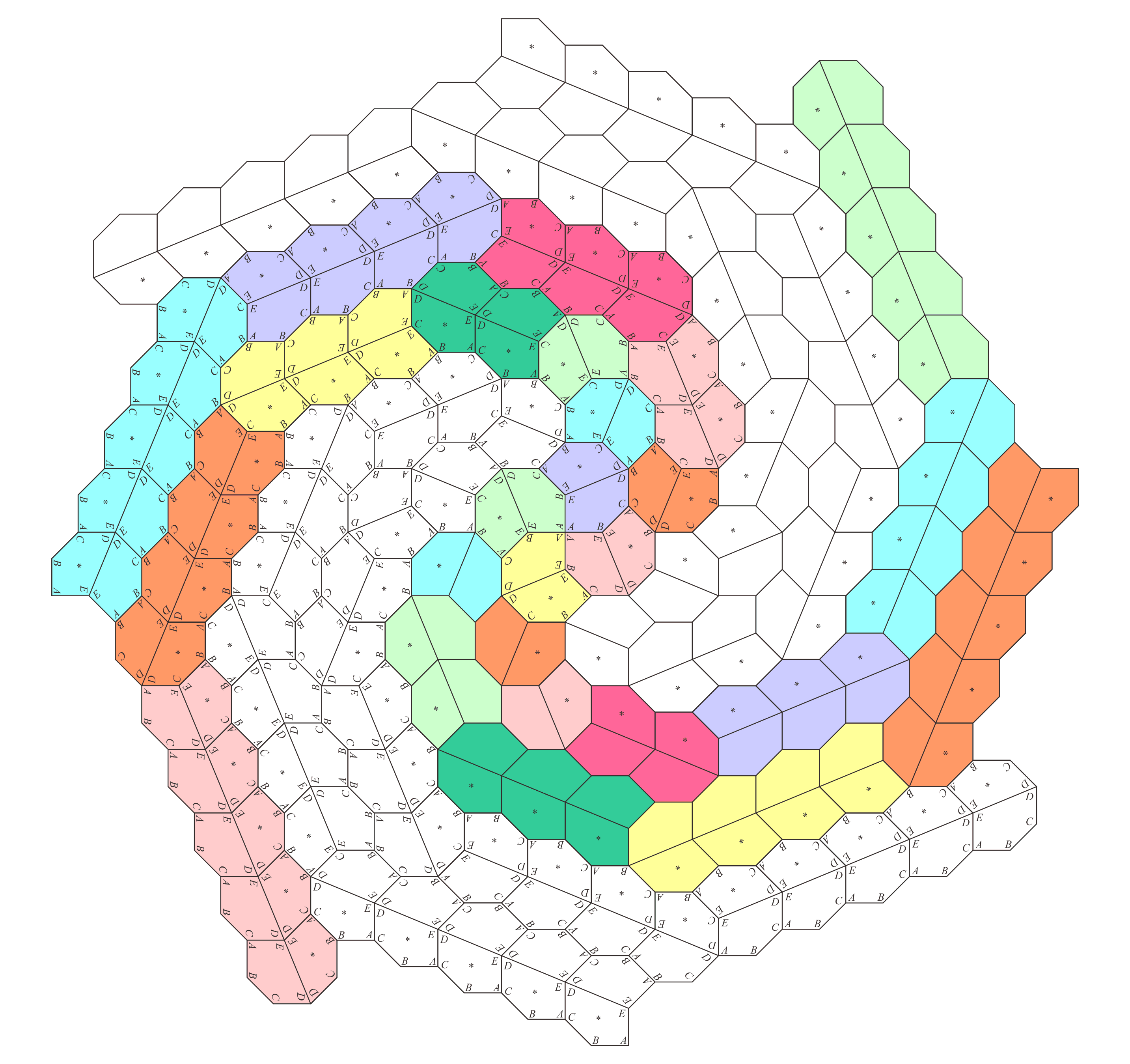} 
  \caption{{\small 
Spiral tiling with four-fold rotational symmetry by a convex pentagonal tile 
\mbox{$P(T1 \cap T7)$} 
} 
\label{fig10}
}
\end{figure}

\begin{figure}[htbp]
 \centering\includegraphics[width=10cm,clip]{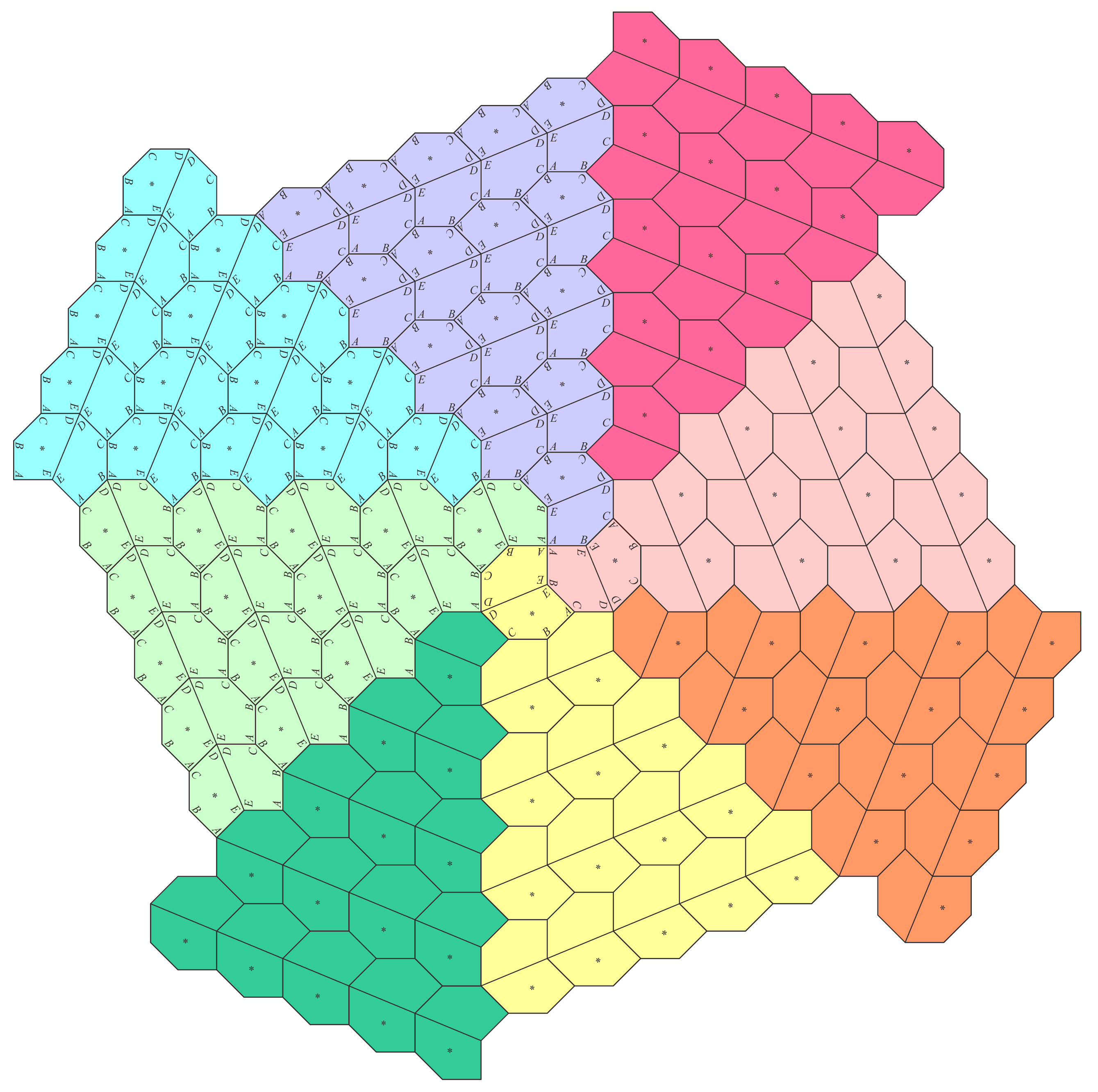} 
  \caption{{\small 
Four-fold rotationally symmetric tiling by a convex pentagonal tile 
\mbox{$P(T1 \cap T7)$}, Part 1
} 
\label{fig11}
}
\end{figure}

\begin{figure}[htbp]
 \centering\includegraphics[width=10cm,clip]{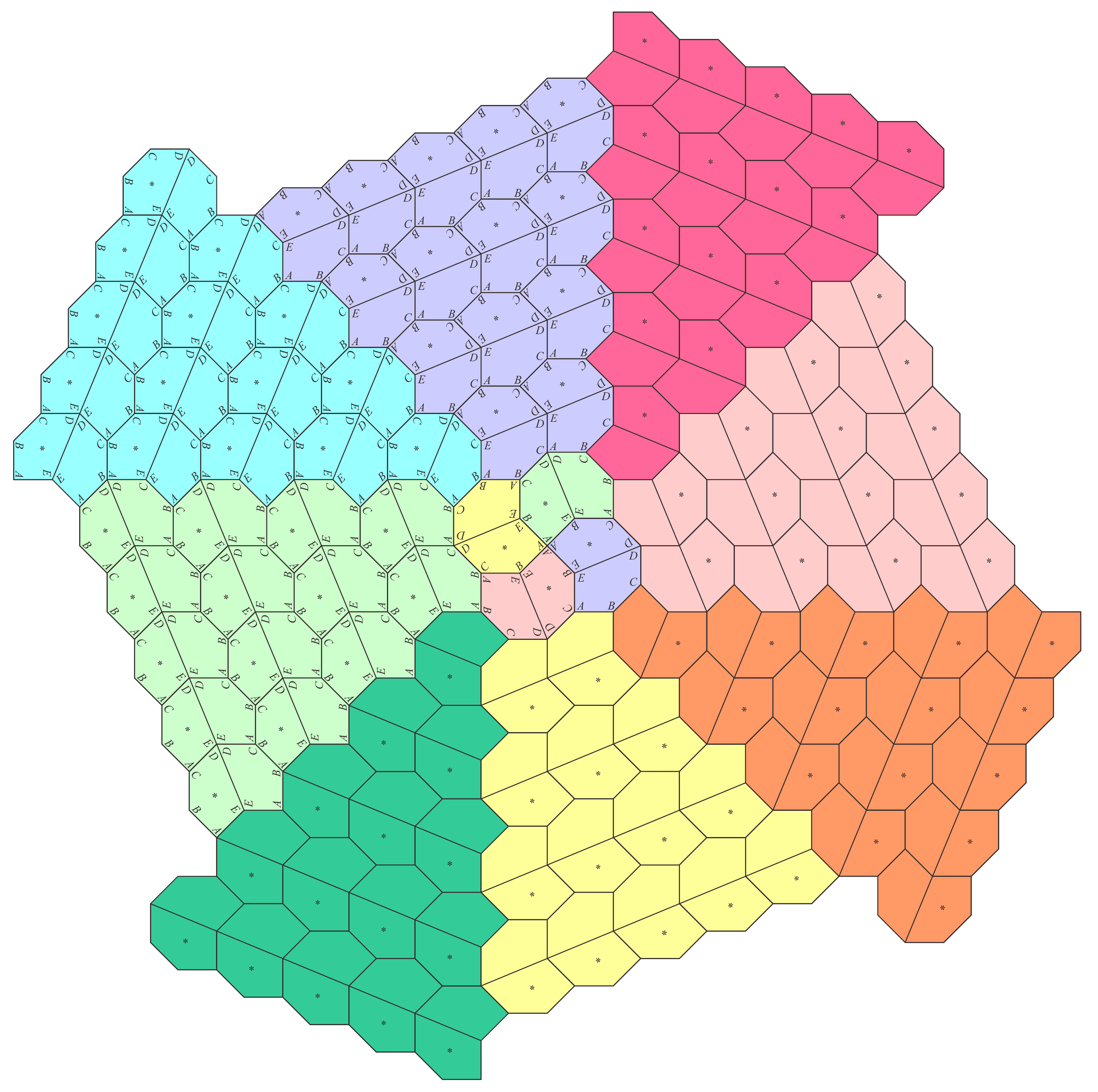} 
  \caption{{\small 
Four-fold rotationally symmetric tiling by a convex pentagonal tile 
\mbox{$P(T1 \cap T7)$}, Part 2
} 
\label{fig12}
}
\end{figure}

\begin{figure}[htbp]
 \centering\includegraphics[width=11.5cm,clip]{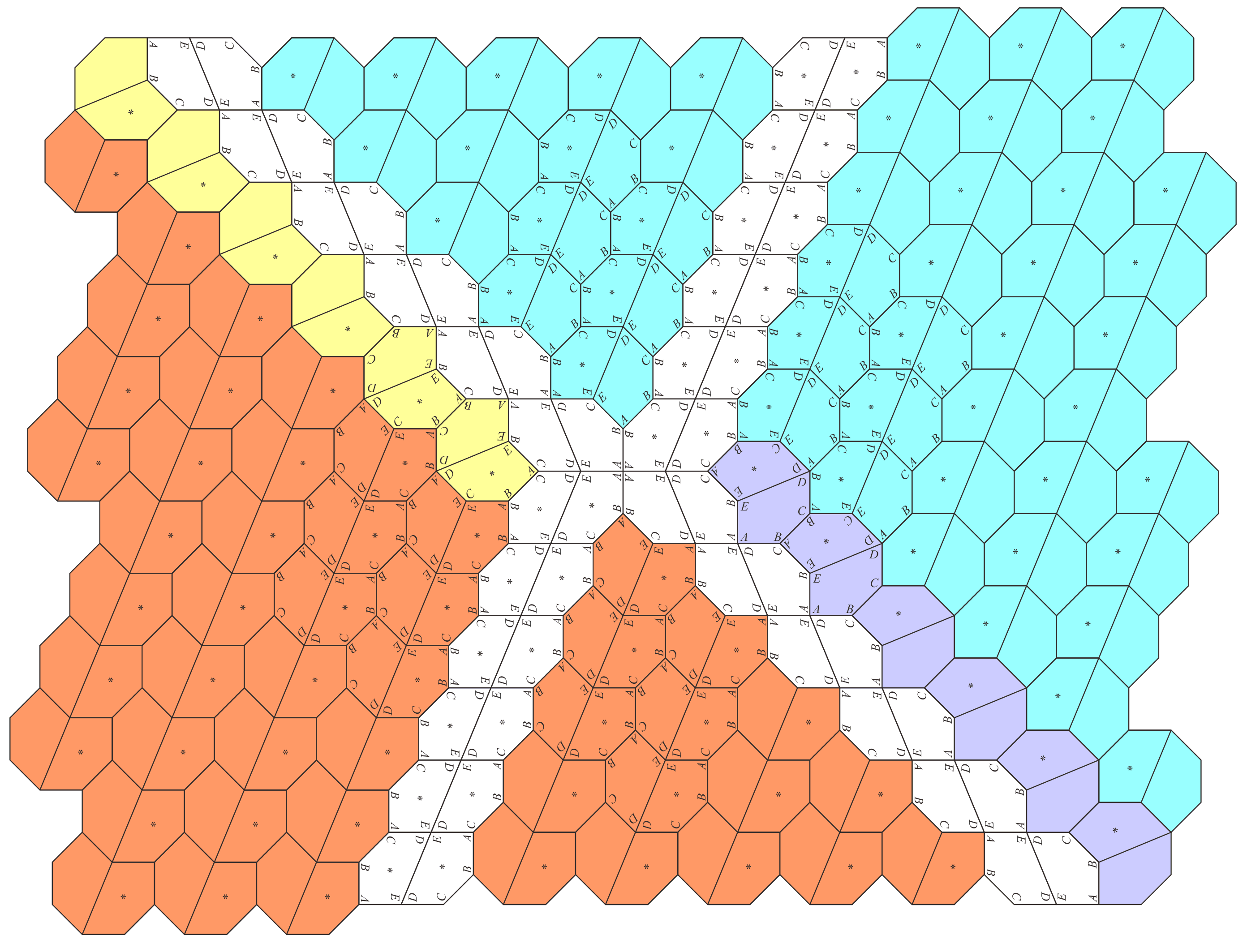} 
  \caption{{\small 
Two-fold rotationally symmetric tiling by a convex pentagonal tile 
\mbox{$P(T1 \cap T7)$}, Part 1 
} 
\label{fig13}
}
\end{figure}

\begin{figure}[htbp]
 \centering\includegraphics[width=11.5cm,clip]{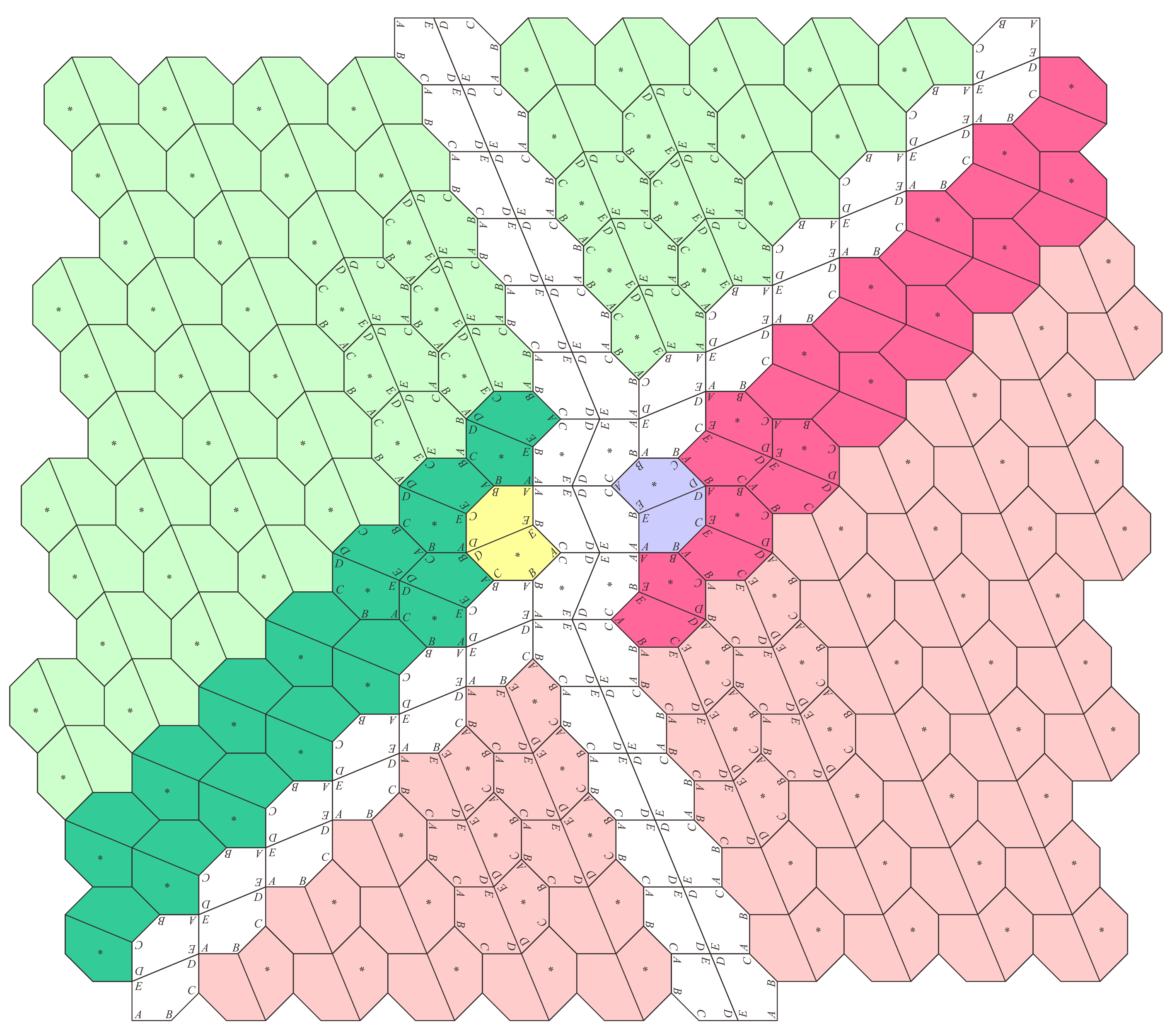} 
  \caption{{\small 
Two-fold rotationally symmetric tiling by a convex pentagonal tile 
\mbox{$P(T1 \cap T7)$}, Part 2
} 
\label{fig14}
}
\end{figure}

\begin{figure}[htbp]
 \centering\includegraphics[width=13.8cm,clip]{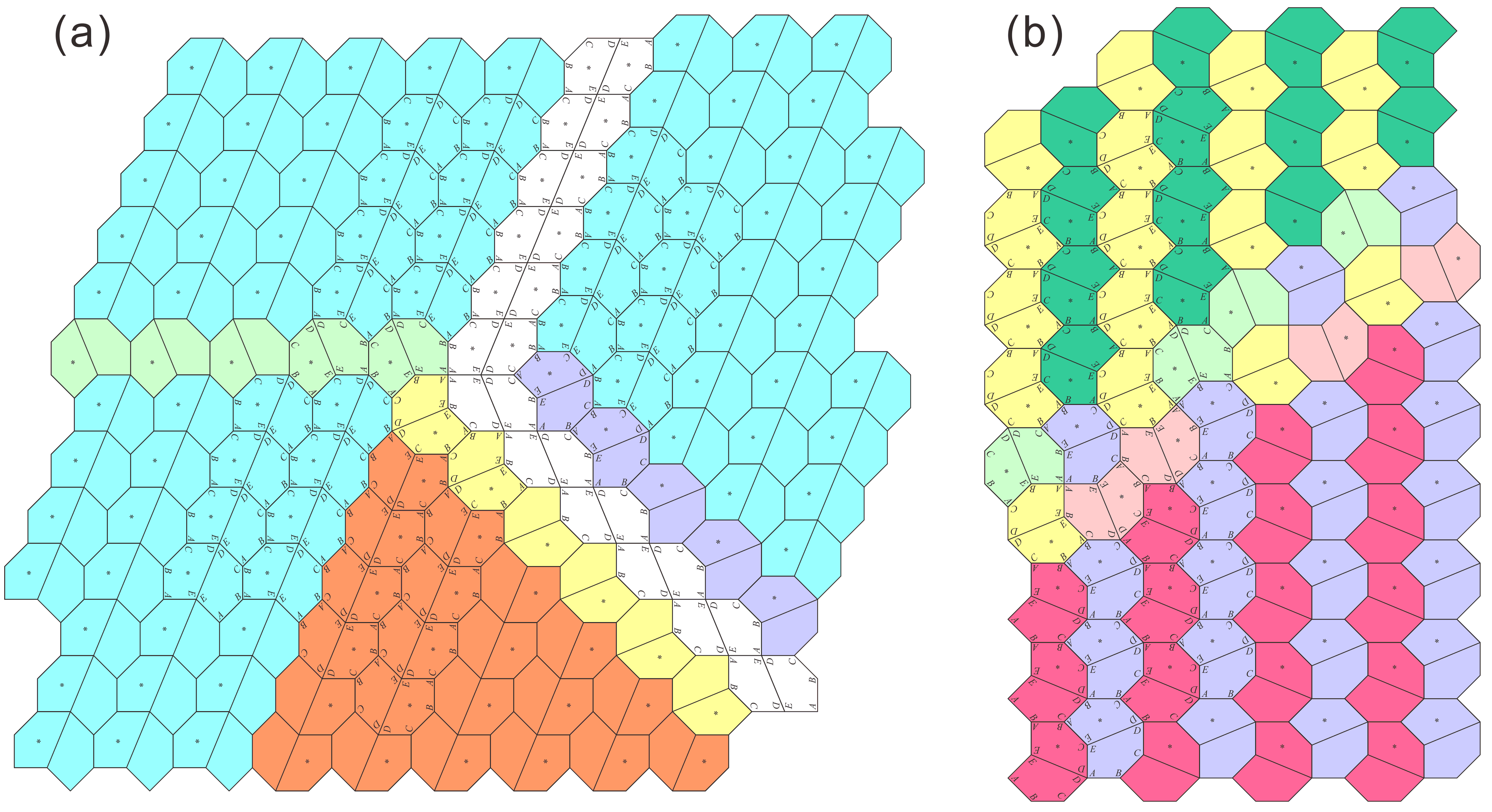} 
  \caption{{\small 
Other tilings with $4A = 360^ \circ $ by a convex pentagonal tile 
\mbox{$P(T1 \cap T7)$} 
} 
\label{fig15}
}
\end{figure}

\begin{figure}[htbp]
 \centering\includegraphics[width=14.2cm,clip]{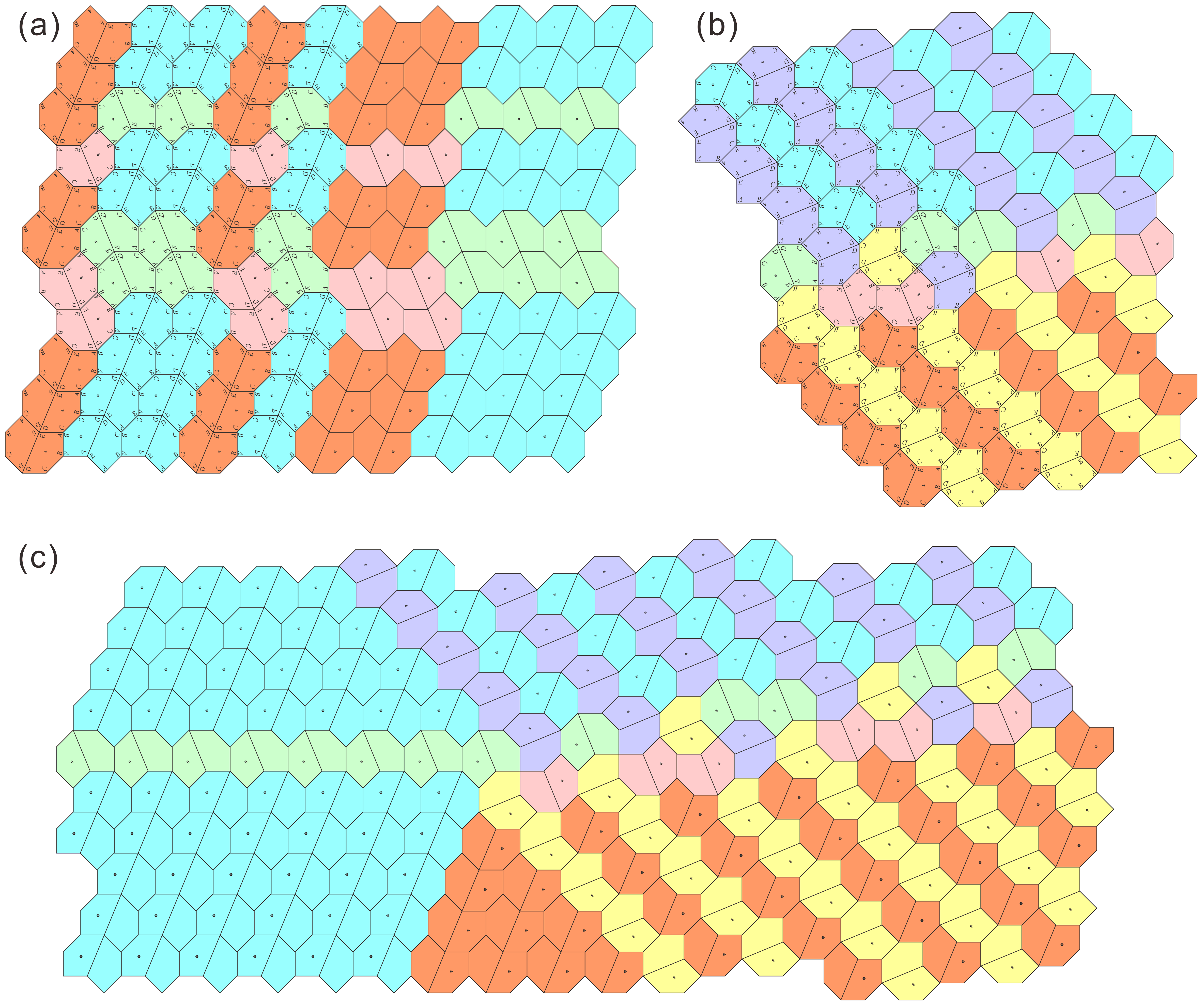} 
  \caption{{\small 
Convex pentagonal tilings with \mbox{$P(T1 \cap T7)$} based on Dailey's other tilings 
} 
\label{fig16}
}
\end{figure}

\mbox{$P(T1 \cap T7)$} can also generate tilings other than those introduced in this 
study. In particular, Yoshiaki ARAKI, of Japan Tessellation Design Association, 
found several interesting tilings using a concave octagon called an Octa-unit. 
They produce four-fold rotationally symmetric tilings that correspond to 
extended tilings such as those shown in Figure~\ref{fig07}, and spiral 
tilings other than the ones described in this study. They will be introduced 
in a separate study.

\bigskip
\noindent
\textbf{Acknowledgments.} 
The author would like to thank Yoshiaki ARAKI of Japan Tessellation Design 
Association, for information and comments.



\renewcommand{\figurename}{{\small Figure.}}

\renewcommand{\figurename}{{\small Figure.}}

\renewcommand{\figurename}{{\small Figure.}}

\renewcommand{\figurename}{{\small Figure.}}

\renewcommand{\figurename}{{\small Figure.}}

\renewcommand{\figurename}{{\small Figure.}}

\renewcommand{\figurename}{{\small Figure.}}

\renewcommand{\figurename}{{\small Figure.}}

\renewcommand{\figurename}{{\small Figure.}}

\renewcommand{\figurename}{{\small Figure.}}

\renewcommand{\figurename}{{\small Figure.}}

\renewcommand{\figurename}{{\small Figure.}}


\begin{thebibliography}{99}
{\small

\bibitem{Dailey_site}
D.~Dailey, A Page for David Daile (\url{http://srufaculty.sru.edu/david.dailey/})
\textgreater tiling (\url{http://cs.sru.edu/~ddailey/tiling/tilingNew.html}) 
\textgreater Dented octagons (\url{http://cs.sru.edu/~ddailey/tiling/octagons.svg}) 
(accessed on 28 March 2020).


\bibitem{G_and_S_1987}
B.~Gr\"{u}nbaum, G.C.~Shephard, \textit{Tilings and Patterns}. 
W. H. Freeman and Company, New York, 1987, pp.15--35 (Chapter 1), 
pp.471--518 (Chapter 9).


\bibitem{Sugimoto_2012}
T.~Sugimoto, Convex pentagons for edge-to-edge tiling, I, 
\textit{Forma}, \textbf{27} (2012) 93--103. Available online: 
\url{https://forma.katachi-jp.com/abstract/2701/27010093.html} 
(accessed on 4 March 2022).


\bibitem{Sugimoto_NoteTP}
\textemdash, Tiling Problem: Convex pentagons for edge-to-edge tiling and 
convex polygons for aperiodic tiling (2015).
\url{https://arxiv.org/abs/1508.01864} (accessed on 23 February 2020).


\bibitem{Sugimoto_2016}
\textemdash, Convex pentagons for edge-to-edge tiling, III, 
\textit{Graphs and Combinatorics}, \textbf{32} (2016) 785--799. 
doi:10.1007/s00373-015-1599-1. 


\bibitem{Sugimoto_2020_1}
\textemdash,  Convex pentagons and concave octagons that can form 
rotationally symmetric tilings (2020).
\url{https://arxiv.org/abs/2005.08470}
(accessed on 19 May 2020).


\bibitem{wiki_pentagon_tiling}
Wikipedia contributors, Pentagonal tiling, 
Wikipedia, The Free Encyclopedia,
\url{https://en.wikipedia.org/wiki/Pentagonal_tiling}
(accessed on 31 March 2020).


\bibitem{wiki_point_group}
\textemdash, Point group, 
Wikipedia, The Free Encyclopedia,
\url{https://en.wikipedia.org/wiki/Point_group}
(accessed on 23 February 2020).


\bibitem{wiki_schoenflies_notation}
\textemdash, Schoenflies notation, 
Wikipedia, The Free Encyclopedia,
\url{https://en.wikipedia.org/wiki/Schoenflies_notation}
(accessed on 23 February 2020).






}

\end{thebibliography}
\end{document}